\documentclass[conference]{IEEEtran}
\IEEEoverridecommandlockouts
\usepackage{amsmath,amssymb,amsfonts,bm}
\usepackage{algorithmic}
\usepackage{graphicx}
\usepackage{mathtools}
\DeclarePairedDelimiter{\ceil}{\lceil}{\rceil}

\begin{document}

\title{A Quantum Annealing Approach for Dynamic Multi-Depot Capacitated Vehicle Routing Problem}


\makeatletter
\newcommand{\linebreakand}{%
  \end{@IEEEauthorhalign}
  \hfill\mbox{}\par
  \mbox{}\hfill\begin{@IEEEauthorhalign}
}
\makeatother

\author{\IEEEauthorblockN{Ramkumar Harikrishnakumar}
\IEEEauthorblockA{\textit{Industrial, Systems, and}\\
\textit{Manufacturing Engineering}\\
\textit{Wichita State University}\\
Wichita, USA \\
rxharikrishnakumar@shockers.wichita.edu}
\and
\IEEEauthorblockN{Saideep Nannapaneni}
\IEEEauthorblockA{\textit{Industrial, Systems, and}\\
\textit{Manufacturing Engineering}\\
\textit{Wichita State University}\\
Wichita, USA \\
saideep.nannapaneni@wichita.edu}
\and
\IEEEauthorblockN{Nam H. Nguyen}
\IEEEauthorblockA{\textit{Boeing Research \& Technology} \\
Huntington Beach, USA \\
nam.h.nguyen5@boeing.com}
\linebreakand

\IEEEauthorblockN{James E. Steck}
\IEEEauthorblockA{\textit{Aerospace Engineering} \\
\textit{Wichita State University}\\
Wichita, USA \\
james.steck@wichita.edu}
\and
\IEEEauthorblockN{Elizabeth C. Behrman}
\IEEEauthorblockA{\textit{Mathematics, 
Physics, and Statistics} \\
\textit{Wichita State University}\\
Wichita, USA \\
elizabeth.behrman@wichita.edu}
}

\maketitle

\begin{abstract}
Quantum annealing (QA) is a quantum computing algorithm that works on the principle of Adiabatic Quantum Computation (AQC), and it has shown significant computational advantages in solving combinatorial optimization problems such as vehicle routing problems (VRP) when compared to classical algorithms. This paper presents a QA approach for solving a variant VRP known as multi-depot capacitated vehicle routing problem (MDCVRP). This is an NP-hard optimization problem with real-world applications in the fields of transportation, logistics, and supply chain management. We consider heterogeneous depots and vehicles with different capacities. Given a set of heterogeneous depots, the number of vehicles in each depot, heterogeneous depot/vehicle capacities, and a set of spatially distributed customer locations, the MDCVRP attempts to identify routes of various vehicles satisfying the capacity constraints such as that all the customers are served. We model MDCVRP as a quadratic unconstrained binary optimization (QUBO) problem, which minimizes the overall distance traveled by all the vehicles across all depots given the capacity constraints. Furthermore, we formulate a QUBO model for dynamic version of MDCVRP known as D-MDCVRP, which involves dynamic rerouting of vehicles to real-time customer requests.  We discuss the problem complexity and a solution approach to solving MDCVRP and D-MDCVRP on quantum annealing hardware from D-Wave. 

\end{abstract}


\section{Introduction}


The vehicle routing problem (VRP) is a well-known combinatorial optimization problem commonly encountered in logistics, transportation, supply chain management, and scheduling \cite{Dantzig1959}. There are several variants of VRP such as multi-depot vehicle routing problem (MDVRP), vehicle routing and scheduling, vehicle routing problem with time windows (VRPTW), capacitated vehicle routing problem (CVRP) and multi-depot capacitated vehicle routing problem (MDCVRP), which are computationally challenging and are NP-hard problems \cite{Milan2017}. In this paper, we discuss a solution approach using the principles of quantum computing to improve the computational performance in solving the MDCVRP and its dynamic version, the D-MDCVRP. 

A MDCVRP can be briefly defined as follows. Given a set of spatially distributed depots with heterogeneous material handling capacities, a variable set of service vehicles at each depot with heterogeneous material carrying capacity, and a set of spatially distributed customer locations, with heterogeneous demands that need to be served, the objective in an MDCVRP is to assign routes to various service vehicles to serve the customer locations considering the depot and vehicle capacities. For example, in a supply chain network for a courier collection service, there will be pickup requests originating from spatially distributed customer locations and these requests need be served by different depots with the help of a heterogeneous vehicle fleet.  The D-MDCVRP is a dynamic version of MDCVRP, which requires dynamic rerouting of vehicles to respond to real-time customer requests.


The Quantum computing paradigm uses the principles of quantum mechanical systems such as superposition and entanglement to improve the algorithmic computational performance when compared to classical counterparts.
Adiabatic quantum computation, introduced by Farhi \cite{Farhi2012}, is based on the the adiabatic theorem \cite{Sarandy2004}. 
Quantum annealing works on the principle of adiabatic quantum computation. Kadowaki-Nishimori \cite{Kadowaki1998} proposed the quantum annealing technique by introducing quantum fluctuations into the classical simulated annealing algorithm for solving combinatorial optimization problems in the traverse Ising model, which is equivalent to a quadratic unconstrained binary optimization (QUBO). Therefore, quantum annealing can be used to solve optimization problems that follow the QUBO formulation. 

Quantum annealing was used to solve a wide variety of combinatorial optimization problems in fields such as bioinformatics \cite{Perdomo-Ortiz2012}, Quantum chemistry \cite{Babbush2014}, computational biology \cite{Li2018}, traffic flow optimization \cite{Neukart2017}, fault diagnosis \cite{Perdomo-Ortiz2015}, training of deep neural networks \cite{Liu2018}, vehicle routing  \cite{Martonak2004}, job shop scheduling\cite{Venturelli2016a} and nurse scheduling in healthcare \cite{Ikeda2019}. Following the application of quantum annealing for solving combinatorial optimization in various domains, we provide the formulations for solving the MDCVRP and D-MDCVRP in a way that facilitates quantum annealing optimization.




 \textbf{\textit{Paper Contributions:}} The contributions made through this paper are: (1) a quadratic unconstrained binary optimization (QUBO) formulation for a MDCVRP under vehicle/depot capacity constraints; (2) a QUBO formulation for dynamic rerouting of multiple vehicles in response to real-time customer requests; and (3) discussion on problem complexity and a solution framework through quantum annealing.

\textbf{\textit{Paper Organization:}} Section \ref{sec:background} provides a background to quantum annealing, QUBO, and D-Wave hardware for quantum annealing. Section \ref{sec:lit_review} provides a review of existing literature on the use of quantum annealing for solving vehicle routing problems. Section \ref{sec:formulation} provides the MDCVRP formulation with constraints along with associated QUBO formulation. Section \ref{sec:dynamic} discusses the D-MDCVRP formulation. Section \ref{sec:discuss} discusses on the problem complexity and solution framework for solving MDCVRP and D-MDCVRP followed by concluding remarks and future work in Section \ref{sec:conclusion}.

\section{Background}
\label{sec:background}
Quantum annealing is a metaheuristic algorithm for solving combinatorial optimization problems based on the principles of quantum mechanics \cite{Kadowaki1998}.  Quantum annealing is implemented by evolving the Hamiltonian dynamics from an initial quantum state to a final quantum state, which corresponds to the solution of an optimization problem of interest \cite{Aharonov2008}. In this context, the dynamism is strictly adiabatic and thus relates to the concept of adiabatic quantum computing \cite{Farhi2001}. If $H_I$ and $H_F$ represent the initial and final Hamiltonians of a quantum system, then the Hamiltonian at any time, $H(t)$, between the evolution can be written as \cite{Farhi2001}:

\begin{equation}
\label{eqn:hamiltonian}
H(t) = \bigg(1-\frac{t}{T}\bigg) H_I + \bigg(\frac{t}{T}\bigg) H_F
\end{equation}

In Eq. \ref{eqn:hamiltonian}, $T$ is the time taken to evolve the system from the initial state to the final state, and $t$ is any time between 0 and T ($0\leq t \leq T$). The quantum system is initially in the ground state of the initial Hamiltonian. When the system is evolved slowly in an adiabatic manner, the system stays in the ground state. The final Hamiltonian, and thus the final ground state, corresponds to the minimization problem of interest.


A quantum system comprises of several individual qubits, and the Hamiltonian of an N-qubit Ising system can be written as:

\begin{equation}
\label{eqn:ising}
H(\bm{s}) = \sum_{1\leq i \leq N} h_i s_i + \sum_{1\leq i \leq j \leq N} J_{ij} s_is_j
\end{equation}

\noindent where $s_i$ represents the state of $i^{th}$ qubit, which can be either $-1$ or $1$, $h_i$ and $J_{ij}$ are the bias and interaction terms relating to individual qubits \cite{su2016quantum}, and $\bm{s}$ is a vector of all qubit states.

In quadratic unconstrained binary optimization (QUBO), the value taken by each individual decision variable is either $0$ or $1$. The Ising model can be translated to a binary optimization representation using the following transformation: $x_i = \frac{1 + s_i}{2}$ by which the $-1$ and $1$ states of $s_i$ are mapped to $0$ and $1$ states of $x_i$. Therefore, Eq. \ref{eqn:ising} can be written as: 

\begin{equation}
\label{eqn:binary}
H(\bm{x}) = \sum_{1\leq i\leq j \leq N} Q_{ij}x_ix_j = \bm{x}^TQ\bm{x}
\end{equation}

In Eq. \ref{eqn:binary}, $Q$ is an $N\times N$ symmetric matrix. Therefore, any optimization problem that can be represented as a QUBO can be solved using quantum annealing. If a constrained binary optimization is available, then it needs to be converted to a QUBO. A commonly used method for such a conversion is the penalty method \cite{Homaifar1994a}. For example, if $f(x)$ and $c(x)\leq 0$ represent the objective function and constraint respectively, then an equivalent unconstrained formulation can be written as $q(x) = f(x)+\gamma g(c(x))$, where $\gamma$ is the penalty term, and $g$ is a function defined over the constraint. In this paper, we discuss QUBO formulations of MDCVRP and D-MDCVRP so they can be solved using quantum annealing.

\section{Literature Review}
\label{sec:lit_review}
 As this paper is considering vehicle routing-related problems, we provide a brief review of previous work that considered quantum annealing approaches to solve such problems.

Boros mentioned a set of search algorithms for solving QUBO problems by illustrating simulation results obtained from various computational experiments \cite{Boros2007}.  According to Choi \cite{Choi2008a}, the adiabatic quantum computation technique can solve QUBO problems that use an Ising-spin Hamiltonian. Furthermore, the adiabatic quantum computation can also solve constrained polynomial optimization problems \cite{Rebentrost2019a}. Vyskocil \cite{Vyskocil2019} proposed how to solve a constrained mixed-integer linear programming problem in the QUBO framework by eliminating large coefficients that often result due to quadratic penalties. The computational techniques such as the adiabatic quantum method, quantum circuits, and quantum walks are suitable to solve Hamiltonian problems but, to address the QUBO problem, we require the application of D-wave architecture platform \cite{Mahasinghe2019a}. 

Moylette \cite{Ikeda2019} proposed a quadratic  speedup quantum algorithm for TSP by using the quantum backtracking algorithm to a classical computation algorithm. The above mentioned previous works show that quantum techniques prove to be beneficial for solving VRP related problems such as MDVRP that motivated towards this work. 

Hirotaka\cite{Irie2019} proposed a QUBO formulation for capacitated vehicle routing problem (CVRP), by introducing the concept of time, capacity and state of vehicles associated with every departure and destinations locations. Sebastian\cite{Feld2019} formulated a hybrid method for solving CVRP using quantum annealer. The heuristic based approach used clustering and routing phase to determine the efficient vehicle routes in each cluster. Clark \cite{Clark2019} investigated providing real-time routing for multiple robots, using hybrid classical-quantum approach for generating collision free routing for multiple robots on simulation grid. Christos\cite{Papalitsas2019} developed a QUBO model formulation for travelling salesman problem (TSP), related to time windows that can handle a small scale TSP with time windows on D-Wave platform.

From the above literature review, we noticed that previous works considered routing problems with a single entity (such as a vehicle or a traveling salesman); however, in real-world applications, we commonly encounter routing problems with heterogeneous vehicles and depots. Therefore, this paper focuses on QUBO formulations for solving MDCVRP and D-MDCVRP on the quantum annealing platform.


\section{MDCVRP Formulation}
\label{sec:formulation}

\subsection{Problem Parameters}
Here, we define the parameters that will later be used in the optimization formulation. 

$x_{ijk}$: A binary variable, which is equal to 1 if location $j$ is served after location $i$ by vehicle $k$.

$\mu_{ik}$: A binary variable, which is equal to 1 if location $i$ is the first location served by vehicle $k$ after leaving its depot.

$\eta_{ik}$: A binary variable, which is equal to 1 if location $i$ is the last location served by vehicle $k$ before returning to its depot.

$\gamma_{kd}$: A binary variable, which is equal to 1 if vehicle $k$ belongs to depot $d$.

$Q_k$: Capacity of vehicle $k$

$q_i$: Demand at location $i$

$V_d$: Capacity of depot $d$

$D_{ij}$: Distance between locations $i$ and $j$

$D_{di}$: Distance distance between depot $d$ and location $i$

$T$: Set of all customer locations

$D$: Set of all depots

$K$: Set of all vehicles across all depots

\subsection{Objective function and Constraints}
\label{subsec:obj_cons}

\textit{Objective function:} We consider minimization of distance traveled by all vehicles across all depots as the objective function. Mathematically, it can be written as 

\begin{equation}
\label{eqn:obj_func}
\begin{aligned}
\text{Min} \hspace{3mm} &\sum_{k \in K} \sum_{i \in T} \sum_{j \in T} D_{ij} x_{ijk} + \sum_{k \in K} \sum_{i \in T} \sum_{d \in D} D_{di} \mu_{ik} \gamma_{kd}\\
& + \sum_{k \in K} \sum_{i \in T} \sum_{d \in D} D_{id} \eta_{ik}\gamma_{kd}
\end{aligned} 
\end{equation}

\textit{Constraint 1:} Each customer location should be served only once across all vehicles. The customer can be the first, somewhere in the middle or at the end of the vehicle route. From a given customer location, a vehicle a visit only one other location.

\begin{equation}
\label{eqn:c1}
\sum_{k \in K} \sum_{\substack{j \in T\\ j \neq i}} x_{ijk} + \sum_{k \in K} \eta_{ik} = 1 \hspace{5mm} \forall i\in T 
\end{equation}

\textit{Constraint 2:} Each customer location can be visited only from one other location across all vehicles from all depots.

\begin{equation}
\label{eqn:c2}
\sum_{k \in K} \sum_{\substack{j \in T\\ j \neq i}} x_{jik} + \sum_{k \in K} \mu_{ik} = 1 \hspace{5mm} \forall i\in T 
\end{equation}

\textit{Constraint 3:} Each trip should start and end at a depot. The variable $\mu_{ik}=1$ indicates that the first customer served by vehicle $k$ after starting from the depot is $i$. Similarly, $\eta_{ik}=1$ indicates that the last customer served by vehicle $k$ before returning to the depot is $i$. Each vehicle starting from a depot should visit only one customer location.

\begin{equation}
\label{eqn:c3}
  \sum_{i \in T} \mu_{ik}=1 \hspace{5mm} \forall k\in K  
\end{equation}

\textit{Constraint 4:} Each vehicle should end at its depot from only a single customer location.

\begin{equation}
\label{eqn:c4}
  \sum_{i \in T} \eta_{ik}=1 \hspace{5mm} \forall k\in K  
\end{equation}

\textit{Constraint 5:} When a vehicle $k$ reaches a location $i$ from $j$, it should leave $i$ to reach some other location $p$. This ensures that a vehicle’s route is continuous and does not terminate at a customer location.
This is often referred to as a flow constraint or continuity constraint.

\begin{equation}
\label{eqn:c5}
\sum_{\substack{j \in T\\ j\neq i}} x_{jik} - \sum_{\substack{p \in T\\ p\neq i}} x_{ipk} = 0  \hspace{5mm}  \forall k\in K, \forall i\in T  
\end{equation}

\textit{Constraint 6:} A vehicle $k$ should not form closed loops with a subset of customer locations. A closed loop (also called a subtour) is formed when the numbers of arcs traveled between any subset of customer locations across all vehicles is equal to the number of customer locations in that subset. Therefore, we constrain the maximum number of arcs to be one less than the number of customer locations. This is known as the subtour elimination constraint \cite{lau2009application}. 

For example, consider a scenario with two depots $D_1$ and $D_2$, and one vehicle at each depot. Let vehicle $V_1$ be associated with $D_1$ and $V_2$ with $D_2$. Assume that there are four customer locations, $C_i, i=1\dots 4$ to be served. When the subtour elimination constraint is not provided, a potential set of routes for vehicles can be as follows. For $V_1$, the route can be $D_1\rightarrow C_1 \rightarrow D_1$ and $C_2 \rightarrow C_3$ and $C_3 \rightarrow C_2$. For $V_2$, the route can be $D_2\rightarrow C_4 \rightarrow D_2$. Here, the vehicles satisfy all the five constraints above, but provide an infeasible route where $V_1$ goes in a loop between $C_2$ and $C_3$. To eliminate such loops, we provide the subtour elimination. This constraint needs to be provided for every subset of customer locations. 

Let $\mathbb{P}(T)$ denote the power set of $T$, i.e., the set of all subsets derived from $T$. Let $S$ represent an element from $\mathbb{P}(T)$, whose cardinality is at least two since formation of a subtour requires atleast two locations. The subtour elimination constraint can be written as

\begin{equation}
\label{eqn:c6}
\sum_{k \in K} \sum_{\substack{i,j \in S \\ i\neq j}} x_{ijk} \leq |S|-1 \hspace{5mm} 2\leq |S|\leq |T|, S \in \mathbb{P}(T)
\end{equation}

\textit{Constraint 7:} Total customer demand should be less than the vehicle capacity

\begin{equation}
\label{eqn:c7}
\sum_{i \in T} \sum_{\substack{j \in T\\ j\neq i}} q_i x_{ijk}  + \sum_{i \in T} q_i\eta_{ik} \leq Q_k \hspace{5mm} \forall k \in K
\end{equation}

\textit{Constraint 8:} Total customer demand across all vehicles should be less than the depot capacity (assuming that the depot capacity is less than the summation of all the vehicle capacities).

\begin{equation}
\label{eqn:c8}
  \sum_{k \in K} \gamma_{kd}\bigg(\sum_{i \in T} \sum_{\substack{j \in T\\ j\neq i}} q_i x_{ijk}  + \sum_{i \in N} q_i\eta_{ik}\bigg) \leq V_d \hspace{5mm} \forall d \in D
\end{equation}

\subsection{Problem Hamiltonian}
\label{subsec:qubo}


Here, we first detail the general approach for the construction of Hamiltonian for any equality and inequality constraint. We later use this approach to construct the Hamiltonian for various constraints in the MDCVRP formulation. First, let us consider an equality constraint of the form, $\sum_{i=1}^{n_x} A_ix_i=b$, where $x_i$ is the $i^{th}$ binary decision variable, $A_i$ represents its coefficient, $b$ is an integer constant, and $n_x$ is the number of decision variables. The Hamiltonan term corresponding to this equality constraint can be written as $(\sum_{i=1}^{n_x} A_ix_i-b)^2$. 

Now, let us consider an inequality constraint of the form, $\sum_{i=1}^{n_x} A_ix_i\leq b$. To represent any optimization formulation as a QUBO, the inequality constraints need to be transformed into equality constraints \cite{Vyskocil2019}; this is accomplished by introducing an additional set of binary decision variables called slack variables, and the corresponding Hamiltonian can be written as $(\sum_{i=1}^{n_x} A_ix_i + \sum_{j=1}^{n_\lambda} 2^j \lambda_j - b)^2$. Here, $n_\lambda$ represents the number of slack variables, which can be calculated as $n_\lambda = \ceil[\big]{1+\log_2 b}$, where $\ceil[\big]{.}$ is the ceiling function. The number of slack variables should be sufficient enough to represent all the values from $0$ (when all the $x_i$ variables are zero) to $b$. Following this discussion, we will write the Hamiltonian terms of various equality and inequality terms in the MDCVRP formulation. We will first start with writing the Hamiltonian term corresponding to the objective function.









The Hamiltonian term, $H_O$ that corresponds to the objective function is simply the objective function in Eq. \ref{eqn:obj_func}.

\begin{equation}
\begin{aligned}
H_O = &\sum_{k \in K} \sum_{i \in T} \sum_{j \in T} D_{ij} x_{ijk} + \sum_{k \in K} \sum_{i \in T} \sum_{d \in D} D_{di} \mu_{ik}\gamma_{kd}\\
& + \sum_{k \in K} \sum_{i \in T} \sum_{d \in D} D_{id} \eta_{ik}\gamma_{kd}
\end{aligned}
\end{equation}


The Hamiltonian terms that correspond to various constraints in Section \ref{sec:formulation} are given below. Let $C_i, i=1\dots 8$ represent the eight constraints, and let $H_{C_i}, i=1\dots 8$ represent the Hamiltonian terms corresponding to the constraints.

\begin{equation}
\label{eqn:hc1}
H_{C_1}=\textit{B}\sum_{i \in T}\bigg (1 - \bigg(\sum_{k \in K} \sum_{j \in T} x_{ijk} + \sum_{k \in K} \eta_{ik} \bigg)\bigg )^2  
\end{equation}

\begin{equation}
\label{eqn:hc2}
H_{C_2}=\textit{B}\sum_{i \in T}\bigg (1 - \bigg(\sum_{k \in K} \sum_{j \in T} x_{jik} + \sum_{k \in K} \mu_{ik} \bigg)\bigg )^2  
\end{equation}


\begin{equation}
\label{eqn:hc3}
H_{C_3}=\textit{B}\sum_{k \in K}\bigg (1 - \bigg(\sum_{i \in T} \mu_{ik} \bigg)\bigg )^2\hspace{5mm} 
\end{equation}


\begin{equation}
\label{eqn:hc4}
H_{C_4}=\textit{B}\sum_{k \in K}\bigg (1 - \bigg(\sum_{i \in T} \eta_{ik})\bigg)\bigg )^2 \hspace{5mm}   
\end{equation}


\begin{equation}
\label{eqn:hc5}
H_{C_5}=\textit{B}\sum_{i \in T} \sum_{k \in K} \bigg(\sum_{\substack{j \in T\\ j\neq i}} x_{jik} - \sum_{\substack{p \in T\\ p\neq i}} x_{ipk}\bigg)^2
\end{equation}


\begin{equation}
\label{eqn:hc6}
\begin{aligned}
H_{C_6} = \sum_{\substack{S \in \mathbb{P}(T)\\ 2\leq |S|\leq |T|}} \bigg( & \sum_{k \in K} \sum_{\substack{i,j \in S\\ i\neq j}} x_{ijk} + \\ 
& \sum_{l=0}^{\ceil[\big]{1+\log_2 |S|-1}} 2^l \lambda_{lS}- |S| + 1 \bigg)^2
\end{aligned}
\end{equation}


\begin{equation}
\label{eqn:hc7}
\begin{aligned}
H_{C_7}=B\sum_{k \in K}\bigg( &\sum_{i \in T} \sum_{\substack{j \in T\\ j\neq i}} q_i x_{ijk}+ \sum_{i \in N} q_i\eta_{ik} + \\
& \sum_{l=0}^{\ceil[\big]{1+\log_2 Q_k}}2^l\lambda_{lk} - Q_k\bigg)^2
\end{aligned}
\end{equation}


\begin{equation}
\label{eqn:hc8}
\begin{aligned}
H_{C_8}=B\sum_{d \in D} \bigg(& \sum_{k \in K} \gamma_{kd}\bigg(\sum_{i \in T} \sum_{\substack{j \in T\\ j\neq i}} q_i x_{ijk}  + \sum_{i \in T} q_i\eta_{ik}\bigg) \\
& + \sum_{l=0}^{\ceil[\big]{1+\log_2 V_d}}2^l\lambda_{ld} - V_d \bigg)^2
\end{aligned}
\end{equation}

In Eqs. \ref{eqn:hc1}-\ref{eqn:hc8},  $B$ is a large positive constant that corresponds to the penalty incurred when the constraints are violated. After obtaining Hamiltonian terms that correspond to the objective function and various constraints, the overall Hamiltonian is equal to the sum of the individual Hamiltonian terms. Therefore, $H_F = H_O + \sum_{i=1}^8 H_{C_i}$. After static formulation, let us now discuss formulation to facilitate dynamic rerouting to respond to real-time customer requests. 

\section{Dynamic MDCVRP formulation for real-time rerouting}
\label{sec:dynamic}

The vehicles across several depots will initially be routed based on the output from the static MDCVRP formulation discussed in Section \ref{sec:formulation}. Assume that at a certain time when are vehicles are in service, a new set of customer requests denoted as $R$ become available. Therefore, the overall set of customer locations, denoted as $W = T\cup R$. Let $\Gamma_W$ and $\Theta_W$ represent the set of locations that were served and yet to be served respectively. Note that $\Gamma_W \subseteq T$, $R \subseteq \Theta_W$, and $W = \Gamma_W + \Theta_W$. When a new set of customer requests arrives and rerouting is performed, the positions of vehicles can be anywhere along the obtained from the static formulation. 

Let $C$ represent the set of current locations of the vehicles. Let $\gamma_{kc}$ represent the binary variable that determines the current location of vehicle $k$.  $\gamma_{kc}=1$ when vehicle $k$ is at location $c$ and $0$ otherwise. These current locations will now become the initial locations of the vehicles for the rerouting process. In the static formulation, the initial locations of various vehicles are their associated depots. Since the vehicles have limited capacity, and some customer locations were already served, the remaining vehicle capacities and similarly, remaining depot capacities needs to be considered in the rerouting analysis. We provide the revised formulation that needs to be adopted for dynamic rerouting. 

The objective function, which is the minimization of total distance to be covered by all the vehicles across all depots can be written as:

\begin{equation}
\label{eqn:d_obj_func}
\begin{aligned}
\text{Min} \hspace{3mm} &\sum_{k \in K} \sum_{i \in \Theta_W} \sum_{j \in \Theta_W} D_{ij} x_{ijk} + \sum_{k \in K} \sum_{i \in \Theta_W} \sum_{c \in C} D_{qi} \mu_{ik}\gamma_{kc}\\
& + \sum_{k \in K} \sum_{i \in \Theta_W} \sum_{d \in D} D_{id} \eta_{ik}\gamma_{kd}
\end{aligned} 
\end{equation}

Constraints 1-6 (Eqs. \ref{eqn:c1}-\ref{eqn:c6}) will have the same form as in the static formulation except for a change in the values taken by various indices. The revised set of constraints are given below.

\begin{equation}
\label{eqn:d_c1}
\sum_{k \in K} \sum_{\substack{j \in \Theta_W\\ j \neq i}} x_{ijk} + \sum_{k \in K} \eta_{ik} = 1 \hspace{5mm} \forall i\in \Theta_W 
\end{equation}

\begin{equation}
\label{eqn:d_c2}
\sum_{k \in K} \sum_{\substack{j \in \Theta_W\\ j \neq i}} x_{jik} + \sum_{k \in K} \beta_{ik} = 1 \hspace{5mm} \forall i\in \Theta_W  
\end{equation}

Here, $\beta_{ik}$ corresponds to the first location reached by a vehicle from the initial position in the rerouting process, and $\beta_{ik}=1$ if location $i$ is the first location served by vehicle $k$.

\begin{equation}
\label{eqn:d_c3}
  \sum_{i \in \Theta_W } \beta_{ik}=1 \hspace{5mm} \forall k\in K  
\end{equation}

\begin{equation}
\label{eqn:d_c4}
  \sum_{i \in \Theta_W} \eta_{ik}=1 \hspace{5mm} \forall k\in K  
\end{equation}

\begin{equation}
\label{eqn:d_c5}
\sum_{\substack{j \in \Theta_W\\ j\neq i}} x_{jik} - \sum_{\substack{p \in \Theta_W\\ p\neq i}} x_{ipk} = 0  \hspace{3mm}  \forall k\in K, \forall i\in \Theta_W  
\end{equation}

\begin{equation}
\label{eqn:d_c6}
\sum_{k \in K} \sum_{\substack{i,j \in S\\ i\neq j}} x_{ijk} \leq |S|-1 \hspace{3mm} 2\leq |S|\leq |\Theta_W|, S \in \mathbb{P}(\Theta_W)
\end{equation}

\begin{equation}
\label{eqn:d_c7}
\begin{aligned}
\sum_{i \in \Theta_W} \sum_{\substack{j \in \Theta_W\\ j\neq i}} & q_i x_{ijk}  + \sum_{i \in \Theta_W} q_i\eta_{ik} \leq Q_k - \\
&\sum_{w \in \Gamma_W} q_w \gamma_{kw} \hspace{5mm} \forall k \in K
\end{aligned}
\end{equation}

In Eq. \ref{eqn:d_c7}, $\gamma_{kw}=1$ if vehicle $k$ served a location $w$ in $\Gamma_W$.

\begin{equation}
\label{eqn:d_c8}
\begin{aligned}
\sum_{k \in K} \gamma_{kd} & \bigg(\sum_{i \in \Theta_W} \sum_{\substack{j \in \Theta_W\\ j\neq i}} q_i x_{ijk}  + \sum_{i \in \Theta_W} q_i\eta_{ik}\bigg) \leq \\
& V_d - \sum_{k \in K}\sum_{w \in \Gamma_W} q_w \gamma_{kw} \gamma_{kd} \hspace{5mm} \forall d \in D
\end{aligned}
\end{equation}

The procedure for writing the Hamiltonian terms corresponding to various objective and constraints in the dynamic rerouting formulation is similar to that used for the static formulation in Section \ref{subsec:qubo}. Since the objective function (Eq. \ref{eqn:d_obj_func}) and the first six constraints (Eqs. \ref{eqn:d_c1}-\ref{eqn:d_c6}) are very similar to the MDCVRP formulation, we will write the Hamiltonian terms for the last two constraints denoted as $H_{C_7}^d$ and $H_{C_8}^d$. Following Section \ref{subsec:qubo}, the number of slack variables required for a given $k$ and $d$ in Eq. \ref{eqn:d_c7} and \ref{eqn:d_c8} can be calculated as $n_{\lambda k} = \ceil[\big]{1+\log_2 (Q_k - \sum_{w \in \Gamma_W} q_w \gamma_{kw})}$ and $n_{\lambda d} = \ceil[\big]{1+\log_2 (V_d - \sum_{k \in K}\sum_{w \in \Gamma_W} q_w \gamma_{kw} \gamma_{kd})}$ respectively. Given $n_{\lambda k}$ and $n_{\lambda d}$, $H_{C_7}^d$ and $H_{C_8}^d$ can be written as  



\begin{equation}
\label{eqn:d_hc7}
\begin{aligned}
H_{C_7}^d=B\sum_{k \in K}\bigg( &\sum_{i \in \Theta_W} \sum_{\substack{j \in \Theta_W\\ j\neq i}} q_i x_{ijk}+ \sum_{i \in \Theta_W} q_i\eta_{ik} + \\
& \sum_{l=0}^{n_{\lambda k}}2^l\lambda_{lk} - Q_k +
\sum_{w \in \Gamma_W} q_w \gamma_{kw}\bigg)^2
\end{aligned}
\end{equation}


\begin{equation}
\label{eqn:d_hc8}
\begin{aligned}
H_{C_8}^d=B\sum_{d \in D} &\bigg( \sum_{k \in K} \gamma_{kd}\bigg(\sum_{i \in \Theta_W} \sum_{\substack{j \in \Theta_W\\ j\neq i}} q_i x_{ijk}  + \sum_{i \in \Theta_W} q_i\eta_{ik}\bigg) \\
& + \sum_{l=0}^{n_{\lambda d}}2^l\lambda_{ld} - V_d +\sum_{k \in K}\sum_{w \in \Gamma_W} q_w \gamma_{kw} \gamma_{kd} \bigg)^2
\end{aligned}
\end{equation}

\section{Discussion}
\label{sec:discuss}
In this section, we discuss the problem complexity and a solution framework for solving the static and dynamic multi-depot capacitated vehicle routing problem. 

\textbf{\textit{Problem Complexity:}} We will use the MDCVRP formulation discussed in Section \ref{sec:formulation} to illustrate the problem complexity. However, the same discussion applied to the dynamic formulation as well. The total number of decision variables is equal to the number of binary decisions that are used to describe various vehicle routes (e.g., $x_{ijk}, \mu_{ik}, \eta_{ik}$) and also the number of slack variables (e.g., $\lambda_{lS}, \lambda_{lk}, \lambda_{ld}$) introduced to convert the inequality constraints to equality constraints. For simplicity, let us refer to the former set of decision variables as route decision variables. The total number of route decision variables is equal to $|T|(|T|-1)|K| + 2|T||K|$, where $|T|(|T|-1)|K|$ is the number of $x_{ijk}$ variables and $|T||K|$ represents the number of each of $\mu_{ik}$ and $\eta_{ik}$.

With regard to slack variables, let us first consider the number of slack variables in the subtour elimination constraint (Eq. \ref{eqn:c6}). We will have a subtour elimination constraint for every $|S|$ combinations of customer locations. Let $\binom{n}{z}$ represent the number of combinations of obtaining $z$ elements from a set of $n$ elements. Therefore, the total number of subtour elimination constraints will be equal to $\sum_{2\leq |S|\leq |T|} \binom{|T|}{|S|}$. For a given value of $|S|$ in the subtour elimination constraint, the number of slack variables is equal to $\ceil[\big]{1+\log_2 (|S|-1)}$. The total number of slack variables considering all the subtour elimination constraints is equal to $\sum_{|S|=2}^{|T|} \binom{|T|}{|S|}\ceil[\big]{1+\log_2 (|S|-1)}$.

Following the discussion on subtour elimination constraint, the number of slack variables in the vehicle capacity constraint (Eq. \ref{eqn:c7}) and depot capacity constraint (Eq. \ref{eqn:c8}) can be calculated as $\sum_{k \in K} {\ceil[\big]{1+\log_2 Q_k}}$ and $\sum_{d \in D} {\ceil[\big]{1+\log_2 V_d}}$ respectively. Thus, the total number of decision variables ($N_D$) is equal to 

\begin{equation}
\label{eqn:nd}
\begin{aligned}
N_D = &|T||K|(|T|+1) + \sum_{2\leq |S|\leq |T|} \binom{|T|}{|S|}\ceil[\big]{1+\log_2 (|S|-1)}\\
&+\sum_{k \in K} {\ceil[\big]{1+\log_2 Q_k}} + \sum_{d \in D} {\ceil[\big]{1+\log_2 V_d}}
\end{aligned}
\end{equation}

\textbf{\textit{Solution framework:}} In order to solve the problem using quantum annealing, the overall Hamiltonian needs to be written in the QUBO form. This can be accomplished using the PyQUBO package in Python \cite{tanahashi2019application}. The PyQUBO compiles the overall Hamiltonian and obtains the coefficient matrix of the binary decision variables, i.e., the $Q$ matrix (in Eq. \ref{eqn:binary}). Depending on the total number of decision variables, the QUBO formulation can be solved either using a quantum solver (such as D-Wave 2000Q) or a hybrid quantum-classical solver \cite{dwavehybrid}. A hybrid solver decomposes the problem into several smaller problems and use a combination of classical and quantum solvers to improve the computational performance. Since D-Wave 2000Q has 2048 qubits, a hybrid solver can handle a higher number of decision variables as it uses a combination of classical and quantum solvers. Also, according to Chancellor \cite{Chancellor2019}, the Chimera graph is not completely utilized due to physical constraints.The analysis steps are summarized below:

\begin{enumerate}
    \item Choose penalty values associated with various constraints; the penalty values may be chosen to be equal for all constraints.
    \item Construct the overall Hamiltonian for the objective and constraint functions by converting any inequality constraints to equality constraints through addition of slack variables.
    \item Obtain the matrix of coefficients (the $Q$ matrix) of the decision variables using the PyQUBO package.
    \item Solve the QUBO problem by choosing an appropriate solver, either a quantum or a hybrid solver, depending on the number of decision variables.
    \item Go to Step 1 if the obtained solution is not desirable.
\end{enumerate}

\section{Conclusion}
\label{sec:conclusion}

This paper discussed quadratic unconstrained binary optimization (QUBO) formulations of Multi-Depot Vehicle Routing Problem (MDVRP) and its variant, the Dynamic MDCVRP (D-MDCVRP) so they can be solved on the quantum annealing hardware such as D-Wave 2000Q. MDCVRP and D-MDCVRP seek to serve a set of spatially distributed customer requests through a heterogeneous vehicle fleet operated from multiple depots with heterogeneous capacities. The MDCVRP and D-MDCVRP  formulations contain both equality and inequality constraints. To derive the QUBO formulations, the inequality constraints need to be converted to equality constraints. In this paper, we introduced an additional of decision variables called the slack variables to convert inequality constraints to equality constraints. We have also discussed the problem complexity in terms of the number of total decision variables and also a step-by-step solution framework to solve the associated QUBO formulations on the quantum annealing hardware. The methods that we discussed in this paper are applicable to combinatorial optimization problems in other domains such as resource management in data centers, warehouses, and microgrids. 


As part of our future work, we will solve the QUBO formulations on the D-Wave 2000Q quantum annealer, and compare the accuracy and computational time of the results against classical heuristic algorithms such as Tabu search. We will also investigate the scalability of the quantum annealing approach by comparing the variation of computational time with respect to the number of vehicles in each depot, and number of depots.




\section*{Acknowledgment}
The authors thank Dr. Ehsan Salari of Wichita State University for helpful discussions regarding the MDCVRP and D-MDCVRP formulations.

\bibliographystyle{IEEEtran}
\bibliography{QCE,OtherRef}

\end{document}